\documentclass[11pt]{amsart}
\usepackage{amsmath,amsthm, amssymb}
  \usepackage[alphabetic]{amsrefs}
\usepackage{graphicx}
\usepackage[backref]{hyperref}

\theoremstyle{plain}
\newtheorem*{main theorem}{Main Theorem}
\newtheorem*{theorem*}{Theorem}
\newtheorem*{emptytheorem*}{}
\newtheorem{theorem}{Theorem}
\newtheorem{lemma}{Lemma}
\newtheorem{sublemma}{Sublemma}[lemma]

\newtheorem*{proposition*}{Proposition}
\newtheorem{corollary}{Corollary}[lemma]

\newtheorem*{conjecture*}{Conjecture}

\newtheorem*{maintechnical}{Main Technical Theorem}
\newtheorem*{minprin}{Minimum Principle}

\theoremstyle{definition}
\newtheorem{definition}{Definition}
\newtheorem*{definition*}{Definition}
\theoremstyle{remark}

\newtheorem*{remark*}{Remark}

\begin{document}
\author{Stefano Luzzatto}
\address{Department of Mathematics,
Imperial College, London}
\email{\href{mailto:stefano.luzzatto@imperial.ac.uk}{stefano.luzzatto@imperial.ac.uk}}
\urladdr{\href{http://www.ma.ic.ac.uk/~luzzatto}
{http://www.ma.ic.ac.uk/\textasciitilde luzzatto}}
\author{Lanyu Wang}
\address{Academy of Mathematics and System Sciences,
Chinese Academy of Sciences}
\email{\href{mailto:lanyu@amss.ac.cn}{lanyu@amss.ac.cn}}

\thanks{The second named author was supported by
China NSF 19901035.}
\date{June 26, 2003}

\title[Topological invariance of generic multimodal maps]{
Topological invariance of generic
non-uniformly expanding multimodal maps}

\begin{abstract}
We show that a strengthened version of the Collet-Eckmann
condition for multimodal maps is topologically invariant. In
particular,  if \( f \) is non-uniformly expanding and the
critical points are generic with respect to the absolutely
continuous invariant measure, then any map \( g \) topologically
conjugate to \( f \) is also
 non-uniformly expanding.
    \end{abstract}
\maketitle

% \tableofcontents
\section{Statement of results}

A main theme of this paper is the relation between the
\emph{topological} structure of a dynamical system and its
\emph{measure-theoretic} and \emph{ergodic} properties. Throughout
the paper we shall be concerned with the class \( \mathcal S \) of
\( C^{3} \) interval maps \( f: I \to I \) with negative
Schwarzian derivative and a finite number of non-flat critical
points with possibly different orders, we give the precise
definitions in section \ref{defs} on page \pageref{defs} below.
Our main result is to show that a sub-class of maps in \( \mathcal S
\)  admit equivalent characterizations
in terms of some topological and metric properties of the critical orbits, 
which in turn imply several ergodic properties such as
the existence of an absolutely continuous invariant measure. This 
subclass is quite abundant in \( \mathcal S \) in the sense that general 
one-parameter families of maps will belong to it for a positive 
Lebesgue measure set of parameters.

\subsection{Non-uniformly expanding maps}
A \emph{physical measure} for \( f \) is any
 ergodic invariant measure \( \mu \) with a basin \( B(\mu \)) of
 positive Lebesgue
 measure, where
 \( B(\mu)=\{x: \lim_{n\to\infty}\frac{1}{n}
 \sum_{i=1}^{n}\delta_{f^{i}(x)} =\mu\} \),
 and  the limit is in the weak-\( \star \) topology.
 Maps in \( \mathcal S \) can admit
 various kinds of physical measures such as Dirac-\( \delta \) measures
 supported on attracting fixed points or absolutely continuous measures.
 There also exist examples for which \( f \) admits no physical
 measures at all \cites{HofKel90b, Kel00}. It was shown in
\cite{BloLyu89} that the number of physical measures
for maps in \( \mathcal S \) is finite and bounded by the
number of critical points.

\begin{definition}
We say that \( f: I \to I \) is
\emph{(non-uniformly) expanding}
if there exists a finite number \( \mu_{1}, .., \mu_{q} \), \( q\geq 1
\), of
ergodic absolutely continuous invariant probability measures
such that the union of their basins has full Lebesgue measure
and such that
\begin{equation}\label{integrability}
0 <  \int \log |Df| d\mu < \infty
\end{equation}
for each \( \mu=\mu_{i} \), \( i=1,.., q \).
\end{definition}

 For simplicity we assume that \( q=1 \) and so there is only one
 absolutely continuous invariant probability measure whose basin has
 full measure and which satisfies
 \eqref{integrability}.
 In the general case our results hold for the restriction of
 \( f \) to the basins \( B(\mu_{i}) \) of each ergodic measure
\( \mu_{i} \).  In fact our results hold for the restriction of \( f
\) to the basin of any ergodic absolutely continuous invariant
probability measure satisfying \eqref{integrability} with no
additional assumptions on the dynamics outside this basin.

\subsection{Topological invariance of non-uniform expansivity}
 We recall that \( f: I \to I \) and \( g: J\to J \)
are \emph{topologically conjugate}, or \( f\sim g \),
if there exists a homeomorphism
\( h: I \to J \) such that \( h \circ f = g \circ h \). A topological
conjugacy preserves topological properties of maps such as dense and
periodic orbits but in principle there is no reason for any
measure-theoretic properties
such as the existence of absolutely continuous
invariant measures to be preserved. We want
to address the following question:

\begin{quote}
  \emph{Suppose that \( f \) is non-uniformly expanding
      and that  \( f \sim g \).  Does this imply
      that \( g \) is non-uniformly expanding ?}
   \end{quote}
A counterexample of Bruin \cite{Bru98a}
shows that the answer is negative in general
even in the unimodal case. However we prove the following
\begin{theorem}\label{generictheorem}
    Suppose that \( f \) is non-uniformly expanding,  \( f \sim g \),
    and that all critical points of \( f
      \) are generic with respect to \( \mu \) in the sense that
      \begin{equation}\label{generic}
\lim_{n\to\infty}\frac{1}{n}\sum_{i=1}^{n} \log |Df_{f^{i}(c)}| =
 \int \log |Df| d\mu.
 \end{equation}
Then \( g \) is non-uniformly expanding.
\end{theorem}

Notice that \eqref{generic} holds for \( \mu \) almost all points and
in fact, by the absolute continuity of \( \mu \) and the assumption
that \( B(\mu) \) has full Lebesgue measure, for Lebesgue almost
all points. This justifies thinking of it as a ``generic'' condition.

We remark that the proof does not imply that the conjugacy \( h \) is absolutely
continuous.  Rather, it is
 based on showing that \( g \) satisfies certain conditions which can
 be shown by independent arguments to imply that \( g \) is
 non-uniformly expanding.

 \subsection{The unimodal case}
The only existing positive results on the question of the topological
invariance of non-uniform expansivity are in the unimodal setting.
 The first example of a (topological class of a )
 map with an absolutely continuous invariant
 measure which is preserved by topological conjugacy was constructed
 by Bruin \cite{Bru94}, another such class is defined in
 \cite{Bru98b}. Sands \cite{San95} formulated a (topological)
 condition on the kneading sequence of the map, i.e. the symbolic itinerary of
 the critical point, and showed that this implied the following
 condition
\begin{equation}\label{nongeneric}
\liminf_{n\to\infty}\frac{1}{n}\sum_{i=1}^{n} \log |Df_{f^{i}(c)}| >
\lambda > 0
\end{equation}
for the (single) critical point \( c \) (and this implies non-uniform
expansivity, see the next section). Notice that \eqref{nongeneric} is
strictly weaker than \eqref{generic}. It
follows from our results that Sands' topological condition is actually
equivalent to a metric condition which we will call Tsujii's condition,
to be defined below, and
which is strictly stronger than \eqref{nongeneric} and strictly
weaker than \eqref{generic}. It was then shown in
\cite{NowPrz98} that \eqref{nongeneric} itself is topologically
invariant, by bringing together results
of previous papers \cites{NowSan98, PrzRoh98}
 and observing  that these
results formed a closed chain of implications which included a condition
which is clearly invariant by topological conjugacy.

These result do not however generalize to the multimodal setting:
a counterexample was recently constructed in \cite{PrzRivSmi03}
exhibiting a pair of quasi-symmetrically (in particular
topologically) conjugate multimodal maps of which one
satisfies \eqref{nongeneric} and the other does not.
We shall show nevertheless that the slightly stronger Tsujii
condition mentioned above is topologically invariant for completely
general multimodal maps. This condition is satisfied by generic maps
\eqref{generic}. First we give a more in-depth discussion
of condition \eqref{nongeneric} and its relation to \eqref{generic}.

\subsection{Collet-Eckmann conditions}
Condition \eqref{nongeneric} is more often formulated as
\begin{equation}\label{exponential}
 |Df^{n}_{f(c)}|\geq C e^{\lambda  n}\quad \text{ for } C,
 \lambda > 0, \ \ \forall\ \ n\geq 1.
\end{equation}
It is traditionally known as the Collet-Eckmann condition since it was first
introduced in \cite{ColEck83} and shown, together with some additional
conditions which were later shown to be unnecessary \cite{NowStr88},
to imply non-uniform expansivity for unimodal maps.
The generalization of these results to the multimodal setting has
been somewhat problematic for various technical reasons and only
recently it has been shown that condition \eqref{exponential} assumed
to hold for \emph{every} critical point does indeed imply non-uniform
expansivity, see
\cite{BruLuzStr03} for the case when all critical points have the same
order and \cite{BruStr01} for the general case.

We remark that in the general multimodal case
condition \eqref{exponential} is always satisfied by Lebesgue
almost every
point if \( f \) is non-uniformly expanding since it is equivalent to
\eqref{nongeneric} which follows from \eqref{generic} which holds
for Lebesgue almost every point by \eqref{integrability} and the
assumptions on \( \mu \).
This justifies the use of the term \emph{expanding} for such a map \( f
\). The \emph{non-uniformity} follows
from the fact that the convergence in \eqref{nongeneric} is
not necessarily uniform and thus in general we have that the constant \(
C \) in \eqref{exponential} is actually a function \( C_{x} \) which
depends on \( x \) and such that
\( \inf_{x}C_{x}=0 \):
it might be necessary to wait for an arbitrary long time
before the exponential behaviour becomes apparent.

The difference between \eqref{generic} and \eqref{nongeneric} is quite
subtle and does not usually play an important role in many argument
which take advantage of the exponential growth property
\eqref{exponential}.
In our results however, this distinction will play a very important
role. We shall define below a set of conditions which together fall strictly
between \eqref{generic} and \eqref{nongeneric}.

\subsection{The Main Theorem}

The statement in Theorem \ref{generictheorem} is somewhat
unsatisfactory due to the asymmetry between the assumptions and the
conclusions: it is not the case in general
that the critical points for \( g \) are generic with respect to the
corresponding absolutely continuous invariant measure for \( g \).
Moreover it does not include some simple, albeit exceptional, cases,
such as when the critical points are non-recurrent, which are clearly
topologically invariant conditions and are
known to imply non-uniform expansivity.
We shall therefore isolate the two main properties of genericity
which are actually used in the proof and show that these properties  are
indeed topologically invariant. One of these properties is the
Collet-Eckmann condition \eqref{nongeneric} (or \eqref{exponential}).
The other is a condition on the recurrence to the critical set.

Let  \( d(x,y)=|y-x| \)  denote the standard Euclidean
distance in \( I \), let \( \mathcal C \) denote the set of critical
points of \( f \) and let
\(d(x)=d(x, \mathcal C) = \min\{d(x,c), \ c\in\mathcal C\}.\)
For any \( \delta > 0
\), let \( \mathcal C_{\delta}=\{x: d(x) \leq \delta\} \)
denote a \emph{metric neighbourhood} of the critical set.
We say that a (critical) point
\( c \) satisfies the \emph{slow recurrence} condition if
\begin{equation}\label{slow}%\tag{SR}
\lim_{\delta \rightarrow 0^{+}} \liminf_{n\rightarrow
+ \infty } \frac{1}{n} \sum_{\substack{1\leq i \leq n\\
f^{i}(c)\in \mathcal C_{\delta}}} \log d(f^{i}(c)) =0.
\end{equation}
Almost every point satisfies  \eqref{slow}
since
\( \log |Df_{f^{i}(x)}| \) and \( \log d(f^{i}(x)) \) are
uniformly comparable by the non-flatness of the critical points
\eqref{nonflat},
Birkhoff's ergodic Theorem implies that
\[
\frac{1}{n} \sum_{f^{i}(x)\in \mathcal C_{\delta}}
\log | Df_{f^i(x)}| \to \int_{\mathcal C_{\delta}}\log |Df| d\mu
\]
and the integrability condition
\eqref{integrability}  implies that
\( \int_{\mathcal C_{\delta}}\log |Df| d\mu \to 0 \) as
\( \delta\to 0 \).
It is therefore in some sense also a generic
condition.

We  say that a map satisfies condition CE if all
its critical points satisfy the Collet-Eckmann condition
\eqref{nongeneric}
and SR if all its critical points satisfy the Slow recurrence
conditions \eqref{slow}. We are now ready to state our main result.

\begin{main theorem}
    Suppose that \( f \) satisfies both CE and SR.
Then every map \( g \) topologically
    conjugate to \( f \) also satisfies both CE and SR.
    \end{main theorem}

Thus the simultaneous occurrence of both conditions CE and SR is
topologically invariant. We stress the fact that neither CE nor SR on
their own are topologically invariant.

\subsection{Tsujii's Conditions}
The measure-theoretical abundance of multimodal maps satisfying
conditions CE and SR simultaneously was proved by Tsujii
\cite{Tsu93}  for quite general families of multimodal maps,
building on a theme of proving the abundance of maps satisfying
various kinds of expansivity and recurrence conditions, all aimed
at ultimately implying some form of non-uniform expansivity
\cites{Jak81, BenCar85, ThiTreYou94, Tsu93, Luz00}. We mention
also that non-uniformly expanding maps are in general very
structurally unstable in the sense that arbitrarily small
perturbations which do not preserve the topological structure of
the map can destroy the expansivity, see \cites{GraSwi97, Lyu02,
Koz03, She, KozSheStr03}.

We conjecture that Tsujii's conditions
might play a similar role in the multimodal setting to that played by
the  Collet-Eckmann conditions in the unimodal case.
Several recent counterexamples show that some key properties of unimodal
CE maps do not extend to the multimodal setting; perhaps these
counterexamples may be overcome by assuming Tsujii's conditions
instead.
Our main theorem is precisely an example in which
this is the case.
Our arguments also indicate possible strategies for proving other
technical properties for maps satisfying Tsujii's conditions,
such as \emph{backward Collet-Eckmann} which is equivalent to CE
in the unimodal case \cite{NowSan98}
but does \emph{not} follow from straightforward
CE in the multimodal case \cite{BruStr01},
\emph{bounded backward contraction}
and \emph{exponential decay of correlations} proved
in \cite{BruStr03} and \cite{BruLuzStr03} respectively for CE
multimodal maps
under the additional assumption that
all critical points have the same order. Counterexamples exist
showing that bounded backward contraction may fail for CE maps with
critical points of different order \cite{BruStr03}. These properties
will be investigated in future papers.

\subsection{Sands' condition}
Our approach is quite distinct from the proof in \cite{NowPrz98}
of the invariance of the CE condition in the unimodal case. The
present strategy is more direct and in the spirit of Sands' thesis
\cite{San95} in which, in the unimodal setting, a topological
condition is formulated in terms of the kneading sequence of the
critical point and this condition is shown to imply the
Collet-Eckmann condition. The relation between Sands' topological
condition and our Slow Recurrence condition \eqref{slow} was
studied in \cite{Wan01} in the unimodal case, and a generalization
of the topological condition given in the bimodal case in
\cite{Wan96}. Here we push this whole approach significantly
further by giving a more general definition of this condition
which in particular applies to the multimodal setting, we call it
a Topological Slow Recurrence (TSR) condition, see \eqref{TSR} on
page \pageref{TSR}. By construction this condition is
topologically invariant and thus our main result follows from
proving the double implication

\begin{maintechnical}
    For any \( f\in\mathcal S \) we have
\begin{equation}\label{CESRTSR}
CE + SR \ \ \Longleftrightarrow \ \ TSR.
\end{equation}
\end{maintechnical}

\subsection{Overview of the paper}

In section \ref{combinatorial} we give all the details of the
combinatorial structure and define condition TSR. We give a brief
summary of the key ideas here.
The critical set \( \mathcal C \) defines a partition \( \mathcal P \) of $I$
which allows us to define a symbolic
itinerary for all points which are not preimages of \( \mathcal C \)
and to define cylinder sets made up of points which share the
same itinerary up to some finite time. The whole structure of cylinder
sets is topological. There is then a natural
\emph{topological} separation time: any point \( x \) near \(
\mathcal C \) shadows the orbit of the critical point for a time
defined by the number of iterates for which \( x \) and \( c \)
have the same symbolic itinerary. This shadowing time tends to \(
\infty \) as \( x \) tends to \( c \). In particular the symbolic
sequence associated to any given point contains information about its
pattern of recurrence to the critical set: if it contains long finite
blocks which coincide with long initial blocks of the symbolic
sequences of one of the critical point then it must have come quite close
to the critical set. This simple idea is used in a crucial way in the
definition of the condition TSR.

In the remaining three sections we proceed to prove the three steps
required to show \eqref{CESRTSR}. We first show that CE + SR implies
TSR, then that TSR implies SR and finally that TSR implies CE.
A particular effort has been made here to make the presentation
as self-contained as possible and  as much as possible independent of
specialized technical notation and constructions
familiar mainly to established specialists in one-dimensional dynamics.

\subsection{S-multimodal maps with non-flat critical points}
\label{defs}
Before ending this section we give the formal definition of the class
\( \mathcal S \) of maps under consideration  and
fix some fairly standard notation and definitions.
We shall  consider  \emph{S-multimodal maps with non-flat critical
points}
by which we mean maps \( f \) satisfying the
 following conditions.
Let $f:I\rightarrow I$ be a $C^3$ map, where $I \subset \mathbb R$
is a closed interval. We shall assume without loss of generality
that \( I=[0,1] \). We call $c \in I$ a critical point of $f$ if
$Df(c)=0$. Denote by $\mathcal C$ the set of all critical points.
$f$ has finite number (say $q$) of \emph{non-flat critical
points} $0<c_1<c_2< \cdots <c_q<1$ if there exist $L_i>1$, $l_i>1$
and a neighbourhood
 \( V(c_{i}) \) of $c_i$ in which
\begin{equation}\label{nonflat}
| x-c_i |^{l_i-1}/L_i \leq |Df(x)| \leq L_i| x-c_i
|^{l_i-1}
\end{equation}
for any $x \in V(c_i)$ and $1 \leq i \leq q$.
$f$ has
\emph{negative Schwarzian derivative}, i.e.
 $ S(f)(x)=({D^3f(x)}/{Df(x)})- ({3 D^2f(x)}/{2 Df(x)})^2
< 0$
for $x \in (0,1) \backslash \mathcal C$. We remark that this
condition implies in particular the non-existence of any inflection
points. Thus all critical points are turning points, i.e. points at
which the maps fails to be a local homeomorphism.

As part of the argument we shall need
to obtain some derivative estimates
about the iterates of \( f \) along the forward orbit of the critical
values. If these orbits map to a periodic point then the estimates
are easily obtained. Otherwise we need to consider
 a two-sided neighbourhood of each critical value and the forward
images of this neighbourhood.
We assume therefore  that
the forward orbits of all the critical
points are disjoint from the boundary of \( I \)
or, if not, that
the boundaries of \( I \) are fixed points.
This is not restrictive since, if it is not the case,
we can extend \( f \)  to a map \( \hat f \) defined
on neighbourhood \( \hat I \) of
\( I \) such that the forward orbits of all the critical
points are disjoint from the boundary of \( \hat I \)
and the dynamics on \( I \) is unaffected.

 \section{Combinatorial structure}
\label{combinatorial}
In this section we describe in  detail the combinatorial
structure associated to any
S-multimodal map with non-flat critical points and having no
neutral or attracting periodic orbits.   We emphasize that
everything we describe here relies on no
other assumptions.
Some of the definitions are closely related to
the so-called \emph{cutting/co-cutting times} used in
\cites{Hof80, HofKel90b, Bru95a, San95, NowSan98} for unimodal maps;
we give a  generalization to the multimodal setting and, we believe,
a simplification of the way those concepts are introduced
and applied.

We fix for the rest of the paper a constant
\begin{equation}\label{delta0}
   \delta_{0}>0
\end{equation}
sufficiently small so that the critical neighbourhood
\( \mathcal C_{\delta_{0}} \) has
precisely \( q \)  disjoint connected components; in particular
any \( x\in \mathcal C_{\delta} \) for \( \delta<\delta_{0} \) is
associated unambiguously to one particular critical point. For any
\( x \) we shall use the notation \( x^{i}=f^{i}(x) \),
especially, $c^i=f^i(c)$ for a critical point $c$.

\subsection{Symbolic dynamics}
The critical points define a natural partition
\[
\mathcal I = \{I_{0}, \dots, I_{q}\}
\]
of \( I \) into \( q+1 \)
subintervals. This allows us to associate
a symbolic \emph{itinerary} to any point
not in the preimage of a critical point
by letting
\[ \underline a(x)
= a_{0}a_{1}a_{2} \ldots
\text{ where }
a_{i}=k  \text{ if }
f^{i}(x) \in I_{k} \text{ for all }  i\geq 0.
\]
For convenience we associate
to each critical point two sequences which differ only in the
first term, so that we think of the partition elements as closed
intervals and think of each critical point as belonging
to both of the adjacent intervals.
We  define the symbolic \emph{separation time} between two
points as
\[
s(x,y) = \min\{k\geq 0: a_{k}(x)\neq a_{k}(y)\}
\]
as long as the itineraries \( \underline a(x) \) and \( \underline
a(y) \) are both defined. If \( y=c \) is a critical
point and \( x \) is in one of the adjacent partition elements
then we naturally consider the  symbolic itinerary for \( c \) for
which \( s(x,y)=s(x,c) \geq 1 \).  In particular we define the time it
takes for a point \( x \) to ``separate'' from the critical set \(
\mathcal C \) by
\[
s(x)= s(x,\mathcal C) = \max\{s(x,c): c\in \mathcal C \}.
\]
Notice that \( s(x) \to \infty \) as \( d(x) \to 0 \).
We  define \emph{topological
neighbourhoods} \( \mathcal C_{n} \) of the critical set \( \mathcal C
\) by
\[
\mathcal C_{n}=\{x: s(x) \geq n\}.
\]
We fix a constant \( N_{0} \) sufficiently large so that
\begin{equation}\label{disjoint}
\mathcal C_{N_{0}}\subset \mathcal C_{\delta_{0}},
    \end{equation}
    recall \eqref{delta0}, and such that
every point \( x \) with \( s(x) \geq N_{0} \) is
associated to some particular critical point \( c \) and
\( d(x)=d(x,c) \leq \delta_{0} \) and \( s(x) = s(x,c) \).

For any given finite sequence \( a_{0}\ldots a_{n} \) with \(
a_{i}\in \{0, \ldots, q\} \)
we define the
\emph{cylinder set}
\[
\hat I^{(n)}_{a_{0} ..a_{n}}=\{x: f^{i}(x) \in I_{a_{i}},
0\leq i \leq n\}
\]
Notice that such a cylinder set may be
empty. For a given point \( x \) with a symbolic itinerary \(
\underline a(x) \) we define the \( n \)'th order cylinder set of \( x
\) as
\[
\hat I^{(n)}(x) = \hat I^{(n)}_{a_{0}\ldots a_{n}}.
\]
This is always a two sided neighbourhood of \(
x \).
We shall often write
\[
\hat I^{(n)}_{-}(x) \quad \text{and}\quad \hat I^{(n)}_{+}(x)
\]
to denote the part of \( \hat I^{(n)}(x) \) to the left and right
respectively of the point \( x \).
Clearly \( \hat I^{(n+1)}(x) \subseteq \hat I^{(n)}(x) \) for any \( n
\). Moreover, by the non-existence of wandering intervals \cites{Guc79, MelStr93}
the preimages of the critical set are dense and therefore
we have
\[
|\hat I^{(n)}(x)| \to 0 \ \text{ as } \ n\to \infty \ \ \forall \ x\in
I. \]
In particular the cylinder sets \( \hat
I^{(n)}(x) \) define a nested sequence of neighbourhoods of the point \(
x \).

\subsection{Partitions}
We now fix some arbitrary point \( x \) not contained in any preimage of the
critical set. We fix some notation related to the partition \(
\mathcal P(x) \) of a neighbourhood of \( x \) defined naturally by
the cylinder sets. More precisely,
consider a
generic cylinder set \( \hat I^{(n-1)}=\hat I^{(n-1)}(x) \).
By definition the images \( f^{j}(\hat I^{(n-1)}) \) all belong to the
same element of \( \mathcal I \)
for all \( j \leq n-1 \). Then there are two possibilities:
\begin{enumerate}
    \item
    \(
f^{n}(\hat I^{(n-1)}) \) does not contain any critical points in its
interior.
In this case all points
continue to share the same combinatorics and we have
\( \hat I^{(n)} = \hat I^{(n-1)} \), and thus in particular
\[ \hat I^{(n-1)}\setminus \hat I^{(n)} =
\emptyset. \]
\item
\( f^{n}(\hat I^{(n-1)})  \) contains one or more critical points
in its interior.
In this case
  the cylinder set \( \hat
I^{(n)} \subset  \hat I^{(n-1)}\) is given precisely by  those
points of \( \hat I^{(n-1)} \) which fall into the same element of \(
\mathcal I \)
as \( f^{n}(x) \) and
\[
\hat I^{(n-1)} \setminus \hat I^{(n)} \neq\emptyset.
\]
In this case \( \hat I^{(n-1)} \setminus \hat I^{(n)} \)
may have  either one or two connected components.
\end{enumerate}
We index the times at which case 2 above occurs
by two sequences
\[ \mathcal N^{-}_{x}= \{n_{i}^{-}\}_{i=1}^{\infty}
\quad\text{ and }\quad
\mathcal N^{+}_{x} = \{n_{i}^{+}\}_{i=1}^{\infty} \]
where
\[ n_{i}^{-} \text{ denotes a time  for which }
\hat I^{(n_{i}^{-}-1)}_{-} \setminus \hat I^{(n_{i}^{-})}_{-}\neq\emptyset
\]
 and
\[
n_{i}^{+}  \text{ denotes a time for which }
\hat I^{(n_{i}^{+}-1)}_{+} \setminus \hat
I^{(n_{i}^{+})}_{+}\neq\emptyset.
\]
For a general point \( x \),
the two sequences
\( \mathcal N^{-}_{x}\) and \(
\mathcal N^{+}_{x} \)
are  independent even though there may well be one or
more pairs \( i, j  \)  such that \( n_{i}^{-}=n_{j}^{+} \).
For the special case of the partition \( \mathcal P(c) \)
associated to one of the critical points
\( c
\), we actually have \( \mathcal N = \mathcal N^{-}_{c}=
\mathcal N^{+}_{c} \).
For each such time \( n_{i}^{\pm} \) we define the intervals
\[
I^{(n_{i}^{\pm})}= \hat I^{(n_{i}^{\pm}-1)}_{\pm} \setminus \hat
I^{(n^{\pm}_{i})}_{\pm}
\]
which define precisely the elements of the partition
\[
\mathcal P(x) = \left\{I^{(n_{i}^{\pm})}: n_{i}^{\pm}\in\mathcal
N_{x}^{\pm}\right\}
\]
of a neighbourhood of \( x \), refer to Figure 1.

\begin{center}
\includegraphics[height=8cm]{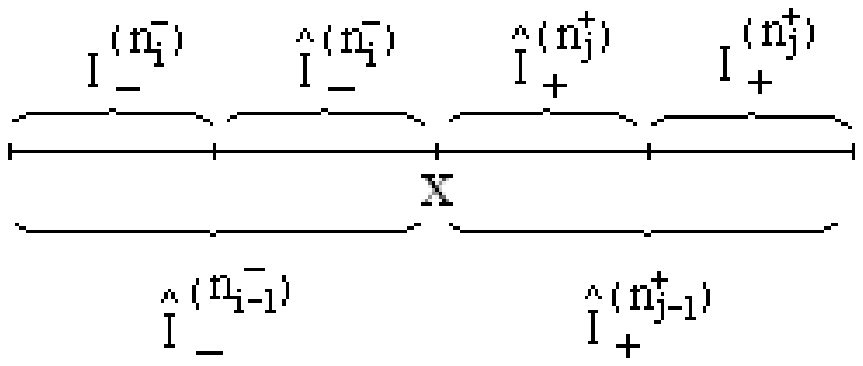}
\makeatletter\def\@captype{figure}\makeatother \caption{Cylinder
sets and partition elements around a point $x$.}
\end{center}

\subsection{Shadowing times}
We call the elements of the sequences \( \mathcal N^{-}_{x}\) and \(
\mathcal N^{+}_{x} \) respectively
the \emph{left and right shadowing times}
 associated to \( x \) because every point in
\( I^{(n_{i}^{\pm})} \) ``shadows'' the point \( x \) for \(
n_{i}^{\pm} \)
iterations:
\begin{equation}\label{separation}
s(y,x) = n_{i}^{\pm} \ \ \Longleftrightarrow\ \ \ y\in
I^{(n_{i}^{\pm})}(x)
\end{equation}
The partition \( \mathcal P(x) \) (and in particular \( \mathcal P(c)
\) for some critical point \( c \)) can be thought of as topological
versions of the partitions defined by the \emph{binding periods}
of Benedicks and Carleson \cite{BenCar85}. Here we lose some good
properties of binding periods such as uniformly bounded distortion
but gain in other ways, as shall become clear below. First we state
three relatively straightforward but not immediately obvious properties
of this partition which will play an important role.
First of all notice that
\begin{equation}\label{monotone}
f^{n_{i+1}^{\pm}} \text{ is monotone on } \hat
I^{(n_{i}^{\pm})}_{\pm}.
\end{equation}
Indeed,
by construction the boundary points of each partition element \(
I^{(n_{i}^\pm)} \) are preimages of some critical point, and one
of the boundary points of each cylinder set $\hat
I_{\pm}^{(n_i^{\pm})}$ is a preimage of some critical point of
order $n_i^{\pm}$. Moreover, there are no points in the interior
of \( \hat I^{(n_{i}^{\pm})}_{\pm} \) which map to a critical
point before time \( n_{i+1}^{\pm} \). Thus \eqref{monotone} follows.

Secondly, for any \( y \) close to \( x \) we have
\begin{equation}\label{distance}
|\hat I^{s(y,x)-1}_{\pm}| \geq d(y,x) \geq |\hat I^{s(y,x)}_{\pm}|.
\end{equation}
Indeed,  by construction we have
\(
\hat I^{(m)}_{\pm}= \hat I^{(n_{i-1}^{\pm})}_{\pm} \) for
\( n_{i-1}^{\pm}\leq m < n_{i}^{\pm}. \)
and so in particular \( \hat I^{(n_{i-1}^{\pm})}_{\pm} =  \hat
I^{(n_{i}^{\pm}-1)}_{\pm}\).
Therefore, for any
\(
y\in
I^{(n_{i}^{\pm})}
\)
the relation \eqref{separation} gives
\[
y\in I^{(n_{i}^{\pm})}=\hat I^{(n_{i-1}^{\pm})}_{\pm} \setminus
\hat I^{(n_{i}^{\pm})}_{\pm} =
\hat I^{(n_{i}^{\pm}-1)}_{\pm} \setminus
\hat I^{(n_{i}^{\pm})}_{\pm}=
\hat I^{s(y,x)-1}_{\pm} \setminus
\hat I^{s(y,x)}_{\pm} \]
which implies  \eqref{distance}.

Thirdly, for  every \( i\geq 1 \) such that \( n_{i+1}^{\pm}
-n_{i}^{\pm} \geq
N_{0} \) we have
 \begin{equation}\label{septime}
s(x^{n_{i}^{
\pm}}) = n_{i+1}^{\pm}-n_{i}^{\pm}.
\end{equation}
To see this consider the two intervals
\(
[y, x] := \hat I^{(n_{i+1}^{\pm})}_{\pm} \subset \hat
 I^{(n_{i}^{\pm})}_{\pm}=: [z, x]
\). By construction we have \( f^{n_{i}^{\pm}}(z) = \tilde c \)
for some critical point \( \tilde c \in\mathcal C \) and thus \(
f^{n_{i}^{\pm}}(\hat I^{(n_{i}^{\pm})}_{\pm})= f^{n_{i}^{\pm}}[z,
x] = [\tilde c, x^{n_{i}^{\pm}}] \) and \( f^{n_{i}^{\pm}}(y) \in
(\tilde c, x^{n_{i}^{\pm}}) \), i.e. the interval \( [\tilde c,
x^{n_{i}^{\pm}}]  \) contains the point \( f^{n_{i}^{\pm}}(y) \)
in its interior. By \eqref{monotone}, \( f^{n_{i+1}^{\pm}} \) is
monotone on \( \hat I^{(n_i^{\pm})}_{\pm} \) which implies in
particular that \( f^{j}(\hat I^{(n_i^{\pm})}_{\pm}) \) cannot
contain any critical point in its interior before time \(
n_{i+1}^{\pm} \), i.e. \( f^{j}(\hat
I^{n_{i}^{\pm}}_{\pm})\cap\mathcal C = \emptyset \) for all \( j <
n_{i+1}^{\pm} \). At this time we have \( f^{n_{i+1}^{\pm}}(y) =
f^{n_{i+1}^{\pm}-n_{i}^{\pm}}(f^{n_{i}^{\pm}}(y)) = \hat c \) for
some (other) critical point \( \hat c \in\mathcal C \). Therefore
$n_{i+1}^{\pm}-n_i^{\pm}$ is exactly the first time at which the
iterates of \( \tilde c \) and of \( x^{n_{i}^{\pm}} \) fall on
different sides of some critical point, and therefore is exactly
the separation time: \( s(x^{n_i^{\pm}}, \tilde c) =
n_{i+1}^{\pm}-n_{i}^{\pm} \). If \( n_{i+1}^{\pm}-n_{i}^{\pm}\geq
N_{0} \)  this implies that \( s(x^{n_i^{\pm}}) = s(x^{n_i^{\pm}},
\tilde c) \) as in \eqref{septime}. See Figure 2.

\begin{center}
\includegraphics[height=8cm]{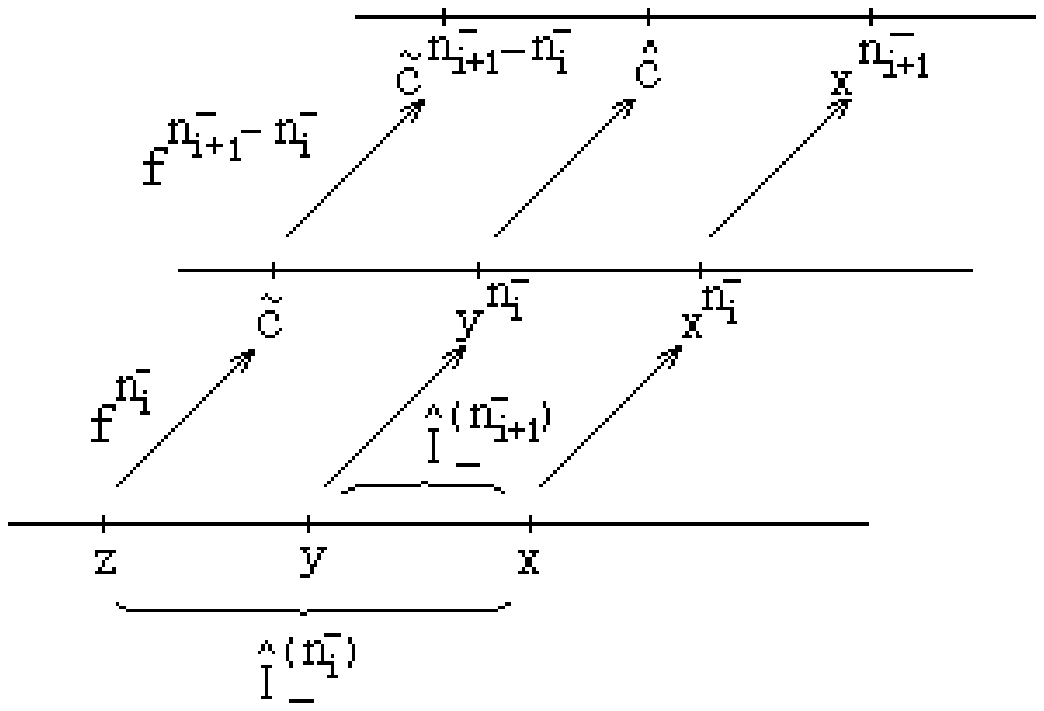}
\makeatletter\def\@captype{figure}\makeatother \caption{Iteration
of the left part of corresponding cylinder set for some left
shadowing time at $x$.}
\end{center}

\subsection{Topological Slow Recurrence} \label{SRTSR}

We are now ready to formulate the \emph{Topological Slow Recurrence}
condition for an arbitrary point \( x \):
\begin{equation}\label{TSR}
    \lim_{m\to\infty} \limsup_{n\rightarrow
+ \infty } \frac{1}{n} \sum_{\substack{1\leq j \leq n\\
x^j \in \mathcal C_{m}}} s(x^{j}) =0.
\end{equation}

\begin{definition}
    We say that the map \( f \) satisfies condition TSR if
    all the critical points satisfy \eqref{TSR}
    \end{definition}
Notice that this depends only on the orbit of the critical points
with respect to the partitions \( \mathcal P(c) \) for \( c\in
\mathcal C \). It is therefore invariant under topological conjugacy.

\section{CE + SR implies TSR}

\begin{lemma}\label{lemma:SRTSR}
    Let \( f \) satisfy condition CE. Then there exists a constant
    \( \bar\kappa>0 \) such that
    \begin{equation}\label{lessthan1}
    s(c^{j}) \leq \bar\kappa \log d(c^{j})^{-1} \quad \forall
    c\in\mathcal C, \ \forall \ j\geq 1.
    \end{equation}
\end{lemma}

We shall use the following
\begin{sublemma}\cite{NowPrz98} \label{NowPrz}
There exists constants \( C, \bar\xi
>0 \) such that for all intervals \( J \) and integers \( s\geq 1 \) such
that \( f^{s}|_{J} \)  is monotone we have
\( |J|\leq Ce^{-s\bar\xi}. \)
\end{sublemma}

\begin{proof}[Proof of Lemma \ref{lemma:SRTSR}]
 We only need to prove the result for \( s(c^j) \geq N_0 \)
 as the choice of the
constant \( \bar\kappa \) can take into account smaller values.
Then there is a well-defined critical point  \( \tilde c \)
``closest'' to \( c^{j} \) i.e. \( s(c^{j}) = s(c^j, \tilde c) \),
and $d(c^j)=d(c^j, \tilde c)=|c^j- \tilde c|$. The map \(
f^{s(c^j)} \) is monotone on the interval \( [c^{j}, \tilde c] \)
and therefore Sublemma \ref{NowPrz} implies \( d(c^{j})=|c^j-
\tilde c|\leq Ce^{-s\bar\xi} \).
\end{proof}

\begin{corollary}\label{cor:TSR}
Let \( f \) satisfy conditions CE and SR. Then it satisfies
condition TSR.
\end{corollary}
    \begin{proof}
    By condition SR, for every \(
 \varepsilon > 0 \) there exist \( n_{\varepsilon},
 \delta_{\varepsilon}>0 \) such that
 \[
 \frac{1}{n} \sum_{\substack{1\leq j \leq n\\
f^{j}(c)\in \mathcal C_{\delta_{\varepsilon}}}} \log
|d(f^j(c)|^{-1} < \varepsilon\quad
\forall \ c\in\mathcal C, \ \ \forall \   n\geq n_{\varepsilon}.
 \]
 Therefore, choosing \( T_{\varepsilon} \) sufficiently large so
 that \( \mathcal C_{T_{\varepsilon}} \subseteq
 \mathcal C_{\delta_{\varepsilon}}  \) and using \eqref{lessthan1}
 we have
 \[
  \frac{1}{n} \sum_{\substack{1\leq j \leq n\\
f^{j}(c)\in \mathcal C_{T_{\varepsilon}}}} s(c^{j}) \leq
  \frac{\bar\kappa}{n} \sum_{\substack{1\leq j \leq n\\
f^{j}(c)\in \mathcal C_{\delta}}} \log |d(c^{j})|^{-1} <
\bar\kappa\varepsilon.
 \]
 Since \( \varepsilon \) is arbitrary, this implies TSR.
   \end{proof}

   \section{TSR implies SR}

\begin{lemma}\label{TSRSR}
    Let \( f \) satisfy condition TSR. Then there exists a constant
    \( \underline\kappa>0 \) such that
    \begin{equation}\label{lessthan}
    s(c^{j}) \geq \underline\kappa\log d(c^{j})^{-1}  \quad \forall
    c\in\mathcal C, \ \forall \ j\geq 1.
    \end{equation}
 \end{lemma}

Arguing as in the proof of  Corollary
\ref{cor:TSR} we then get
\begin{corollary}
 Let \( f \) satisfy condition TSR.
 Then it satisfies condition SR.
    \end{corollary}

We first give two preliminary results which will be used in the proof.
 \begin{sublemma}\label{exp}
     There exists constants
     \( C, \underline\xi >0 \) such that for any interval \( J \) and
     any positive integer \( n \)  such that \( J \) is a
     one-side neighbrouhood of a critical point \( c \),  \(
     f^{n}(J) \)  contains a one-sided
     neighbourhood of the same critical point \( c \)
     and  \( f^{n+1}(J)\supset f(J) \),
     we have \(     |J| \geq C e^{-n\underline\xi}
     \) (we do not assume that \( f^{n+1}|_{J} \) is monotone).
     \end{sublemma}
     \begin{proof}
By the non flatness of the critical
points, there exists a constant \( \theta > 0 \) depending only on \(
f \) such that in either case we have \( \theta |J| \leq |f^{n}(J)|  \).
Moreover, letting \( D=\max_{x\in I}|Df(x)| \),  the
mean value theorem implies \( |f^{n}(J)| \leq D^{n-1} |f(J)|\), and,
using the non-flatness of the critical points again, we have \( |f(J)| \leq
L |J|^{\ell}\) for some constants \( \ell, L \) depending only on
\( f \). Combining these estimates give  \( \theta |J|\leq
 L D^{n-1}|J|^{\ell}  \) and therefore \(
|J|^{\ell-1}\geq \theta L D^{-n +1} \) or \(  |J| \geq (\theta
L^{-1}D)^{1/\ell-1}
D^{-n/(\ell-1)}\).
\end{proof}

\begin{sublemma}\label{integer}
   Let $f$ satisfy condition TSR. Then for any \( \varepsilon > 0 \)
   there exists \( N_{\varepsilon}>0
\) such that
\begin{equation}\label{epsilon}
n_{i+1}(c) -n_{i}(c) \leq \varepsilon n_{i}(c) \quad \forall \
 c\in\mathcal C , \ \ \forall \
n_{i}\geq N_{\varepsilon}.
\end{equation}

\end{sublemma}
    \begin{proof}
From TSR we have \( s(c^{j})/j \to 0\)
as \(j\to \infty \). By \eqref{septime} this immediately implies the
statement.
\end{proof}

\begin{proof}[Proof of Lemma \ref{TSRSR}]
We first define in an iterative process which may stop according to one of
two possible stopping rules. We then show
that in both cases the conclusions of the Lemma are
satisfied. We fix some \( \varepsilon\in (0,1) \) and let \(
N_{\varepsilon} \) be the corresponding integer from Sublemma
\ref{integer}.
\subsection*{Step 1}
Let
\(
s_{0}= s(c^{j}).
\)

\textbullet\emph{ If \( s_0 <
 \max \{N_{0}, N_{\varepsilon} \} \) go to step 3.}

 \noindent
 Otherwise we argue as follows.
Since \( s_{0}\geq N_{0} \), there is a well
 defined critical point \( c_{(0)}
 \) such that \( d(c^{j}) = d(c^{j}, c_{(0)}) \)
 and \( s(c^{j}) = s(c^{j}, c_{(0)}) \) . By \eqref{separation}
 there exists some \( n_{i_{0}}\in\mathcal N_{c_{(0)}} \) with
 \(
 n_{i_{0}}=s_{0}
 \) and
 \[ \hat I^{n_{i_{0}}} \subset (c^{j}, c_{(0)}) \subset \hat I^{n_{i_{0}-1}}. \]
%
%  so by \eqref{distance} we have that
%  \[
%  |\hat I^{n_{i_{0}-1}}|\geq d(c^{j}) \geq |\hat I^{n_{i_{0}}}| >
%  |\hat I^{n_{i_{0}+1}}|
%  \]
 Here we omit the subscripts \( \pm \) not to overload the notation,
 and let $\hat I^{n_{i_0}}=\hat I_{+}^{n_{i_0}}$ or $\hat I_{-}^{n_{i_0}}$
according to the relative positions of \( c^{j} \) and \( c_{(0)} \).
%  Thus it is sufficient to prove
%  \[
%  |\hat I^{n_{i_0}}| = |\hat I^{s_{0}}|\geq e^{-\underline\kappa s_{0}}.
%  \]
 By construction
 \[
 f^{s_{0}}( \hat I^{n_{i_0}}) = (c_{(0)}^{s_{0}}, c_{(1)})
 \]
for some critical point \( c_{(1)} \).

\textbullet\emph{If \( c_{(1)}=c_{(0)} \)  go to step 2.}

\noindent
Otherwise we repeat the algorithm
 with \(  c_{(0)}^{s_{0}}\) playing exactly the role of \( c^{j} \)
 above. More precisely, we consider the
 separation time
 \(
 s_{1}=s(c_{(0)}^{s_{0}})
%  = s( c_{(0)}^{s_{0}}, c_{(1)}).
 \)
 Since \( s_{0}>N_{\varepsilon} \),  \eqref{septime} and \eqref{epsilon} give
\[
 s_{1}= n_{i_{0}+1}-n_{i_{0}} \leq \varepsilon n_{i_{0}} =
 \varepsilon s_{0}.
\]

\textbullet\emph{ If \( s_{1}< \max\{N_{0}, N_{\varepsilon}\} \)  go to step 3.}

\noindent
Otherwise there exists some
 \(
 n_{i_{1}} \in\mathcal
 N_{c_{(1)}} \) such that
 \(  n_{i_{1}}=s_{1} \) and
 \[
 \hat I^{n_{i_{1}}} \subset (c_{(0)}^{s_{0}}, c_{(1)}) \subset
 \hat I^{n_{i_{1}-1}}.
 \]
Again, by construction
 \[
 f^{s_{1}}( \hat I^{n_{i_{1}}}) = (c_{(1)}^{s_{1}}, c_{(2)})
 \]
for some critical point \( c_{(2)} \).

\textbullet\emph{If \( c_{(2)}\) equals either \( c_{(0)} \) or \( c_{(1)} \)
 go to step 2. }

 \noindent
 Otherwise we repeat the process
 again with \( c^{s_{(1)}}_{(1)} \) playing the role of \( c^{s_{(0)}}_{(0)}
 \) and continue in this way until we go to either step 2 or step 3.
 In the general case we have \( s_{j}=s(c_{(j-1)}^{s_{j-1}}) \)
 and
 \begin{equation}\label{epsilon1}
     s_{j}=n_{i_{j-1}+1}-n_{i_{j-1}} \leq \varepsilon n_{i_{j-1}} \leq
 \varepsilon^{j}s_{0}.
 \end{equation}

 \textbullet\emph{If \( s_{j}< \max\{N_{0}, N_{\varepsilon}\} \)  go to step 3.}

 \noindent Otherwise there exists some \(
 n_{i_{j}}\in\mathcal N_{c_{(j)}} \) with \( n_{i_{j}}=s_{j} \) such that
 \[
  \hat I^{n_{i_{j}}} \subset (c_{(j-1)}^{s_{j-1}}, c_{(j)}) \subset
 \hat I^{n_{i_{j}-1}}
\text{  and }
 f^{s_{j}}( \hat I^{n_{i_{j}}}) = (c_{(j)}^{s_{j}}, c_{(j+1)})
 \]
for some critical point \( c_{(j+1)} \).

\textbullet\emph{If \( c_{(j+1)}\) equals any one of
\( c_{(0)}, c_{(1)}, \ldots, c_{(j)} \)
 go to step 2. }

 \noindent The process has to stop in a maximum of \( q+1 \) steps
 as by that time the above condition is necessarily satisfied. See
 Figure 3.

\begin{center}
\includegraphics[height=8cm]{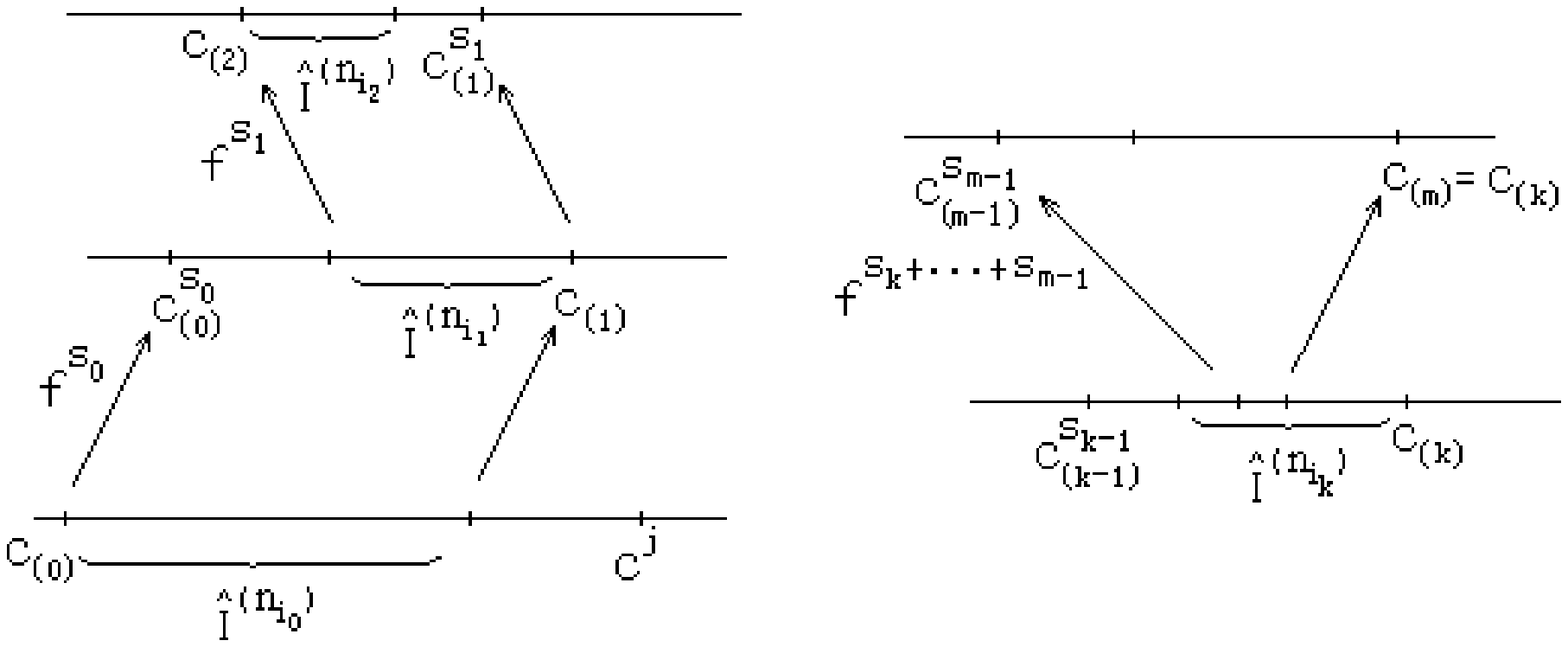}
\makeatletter\def\@captype{figure}\makeatother \caption{Iterative
process at the critical point, (a) step 1, (b) step 2.}
\end{center}

 \subsection*{Step 2}
We suppose here that the procedure described above gives
constants   \(
 0\leq k < m \leq q+1\)
 such that \( c_{(m)}=c_{(k)} \) and therefore we have one-sided
 neighbourhoods
 \( (c_{(k-1)}^{s_{k-1}}, c_{(k)}) \) and
 \( (c_{(m-1)}^{s_{m-1}}, c_{(m)}) =
 (c_{(m-1)}^{s_{m-1}}, c_{(k)}) \) of the same critical point \( c_{(k)}
 \).
 Moreover, by  construction we have that
 \( f^{s_{k}+\ldots+s_{m-1}}(c_{(k-1)}^{s_{k-1}}, c_{(k)})
 \supset
 (c_{(m-1)}^{s_{m-1}}, c_{(m)}) =
 (c_{(m-1)}^{s_{m-1}}, c_{(k)})\) and therefore we are in a position
 to apply the estimates of Sublemma \ref{exp} to the interval
 \( J= (c_{(k-1)}^{s_{k-1}}, c_{(k)})  \)
 if we can show that
  \begin{equation}\label{sep}
      f(c_{(m-1)}^{s_{m-1}}, c_{(k)})
  \supset f(c_{(k-1)}^{s_{k-1}}, c_{(k)}).
  \end{equation}
 To see that this is the case we recall that
 \( s_{j} \) satisfies \eqref{epsilon1} for all \( 0\leq j < m
 \) (for otherwise we would have gone straight to Step 3) and thus in
 particular  the separation time \( s(c_{(m-1)}^{s_{m-1}}, c_{(k)})  =
 s(c_{(m-1)}^{s_{m-1}}, c_{(m)})= s_{m}\) is
 strictly less than \(  s(c_{(k-1)}^{s_{k-1}}, c_{(k)}) =s_{k}\).
 This implies that the separation time \(  s(f(c_{(m-1)}^{s_{m-1}}), f(c_{(k)}))\)
 is strictly less than \(  s(f(c_{(k-1)}^{s_{k-1}}), f(c_{(k)})) \)
 where both \( (f(c_{(m-1)}^{s_{m-1}}), f(c_{(k)})) \) and
 \( (f(c_{(k-1)}^{s_{k-1}}), f(c_{(k)})) \) are one-sided neighbourhood
 (on the same side) of the critical value \( f(c_{(k)}) \). Clearly this
 implies that \( |f(c_{(k-1)}^{s_{k-1}})- f(c_{(k)})| <
 |f(c_{(m-1)}^{s_{m-1}}) - f(c_{(k)})|\) which is precisely \eqref{sep}.
%  In particular
%  \begin{equation}\label{INC}
% \hat I^{(n_{i_{k}})}_{\star}
%   \subset
%   (c_{(k-1)}^{s_{k-1}}, c_{(k)}) \subset
%    \hat I^{n_{i_{k}-1}}
%    \subset (c_{(m-1)}^{s_{m-1}}, c_{(m)}) =
%    (c_{(m-1)}^{s_{m-1}}, c_{(k)}).
%  \end{equation}
% $\star =+$ or $-$. Notice that $|\hat I^{(n_{i_{k}})}_{+}|$ is
% comparable to $|\hat I^{(n_{i_{k}})}_{-}|$ by the non-flatness of
% critical points, there exists a constant $D_1>1$ such that $|\hat
% I^{(n_{i_{k}})}_{+}|/D_1 \leq |\hat I^{(n_{i_{k}})}_{-}|\leq
% D_1|\hat I^{(n_{i_{k}})}_{+}|$. By \eqref{INC} we apply sublemma
% \ref{exp} to $\hat I^{n_{i_k}}=\hat I^{(n_{i_k})}_{\star}$, we
% have
Therefore, by Sublemma \ref{exp} we have
 \begin{equation}\label{1}
    |(c_{(k-1)}^{s_{k-1}}, c_{(k)})| \geq C e^{-\underline\xi (s_{k}+ \ldots + s_{m-1})}.
 \end{equation}
 By construction we also have \( f^{s_{0}+\ldots+s_{k-1}}(c^{j},
 c_{(0)}) \supset  (c_{(k-1)}^{s_{k-1}}, c_{(k)}) \) and therefore by
 the mean value theorem we have
 \begin{equation}\label{2}
 |(c^{j},  c_{(0)}) | \geq D^{-(s_{0}+\ldots+s_{k-1})}
 |(c_{(k-1)}^{s_{k-1}}, c_{(k)})|.
 \end{equation}
 Combining \eqref{1} and \eqref{2} we get that there exist constants
 \( C, \eta > 0 \) such that
 \(
  d(c^{j}) = |(c^{j},  c_{(0)})| \geq Ce^{-\eta (s_{0}+\ldots +
  s_{m-1})}.
 \)
 Finally, from \eqref{epsilon1} we have that there exists
 a uniform constant \( \nu = 1+\sum \varepsilon ^{i} \) such that
\( s_{0}+ \ldots + s_{m-1} \leq \nu s_{0}= \nu s(c^{j}), \) which
then gives \(
  d(c^{j}) = |(c^{j},  c_{(0)})| \geq Ce^{-\eta \nu s(c^{j})}.
 \) This completes the proof in this case.

\subsection*{Step 3}
We now consider the situation in which the procedure described in
Step 1 leads to the existence of some
$0 \leq j \leq q+1$
 such that $s_j< \max \{ N_{0}, N_{\varepsilon} \}$.
Since
 $s_j=s(c^{s_{j-1}}_{(j-1)}, c_{(j)})$ we have that $|c^{s_{j-1}}_{(j-1)}-
  c_{(j)}|\geq \delta$ where \( \delta>0 \) is some constant
  depending only on \( f \), \( N_{0} \) and \( N_{\varepsilon} \).
Therefore a similar argument to that used in the final part of Step 2
gives the result in this case also.
\end{proof}

\section{Topological Slow Recurrence implies CE}
\label{TSRCE}

We now show that for
S-multimodal maps the Topological Slow
Recurrence condition implies the Collet-Eckmann condition.
We shall be concentrating here on
the cylinder sets associated to the
\emph{critical values} rather than those associated to the critical
points used  in the arguments given in the previous section.
 Consider a critical value $c^1$ with shadowing times $n_i^{\pm}$
and corresponding intervals
 $\hat I_{\pm}^{(n_i^{\pm})}$ and  $I^{(n_i^{\pm})}$
 as in section \ref{combinatorial}.

\begin{definition}
For any (large) \( T > 0 \) we say that the gap
 \( n_{i}^{-}-n_{i-1}^{-} \), respectively \( n_{j}^{+}-n_{j-1}^{+} \),
between two left, respectively  right, shadowing times is
\emph{short} if \( n_{i}^{-}-n_{i-1}^{-}< T \), respectively \(
n_{j}^{+}-n_{j-1}^{+}< T \). We say that a left short gap and a
right short gap are \emph{simultaneous} if they are contained in
an interval \( \mathcal T = [t_{1}, t_{2}] \) of times of length
\( t_2-t_1 < T \).
\end{definition}

In section \ref{density} we show that there exist values of \( T
\) for which there are many intervals of simultaneous short
shadowing times. In section \ref{step1} we show that each such
interval implies a uniform estimate on the exponential shrinking
of the length of intervals associated to some particular cylinder
sets at the corresponding critical value. Finally, in section
\ref{global} we combine these two results to obtain the
Collet-Eckmann property.

\subsection{Positive density of simultaneous short shadowing times}
\label{density}

\begin{lemma}\label{posden}
For every small \( \varepsilon>0 \) there exists $T=T_{\varepsilon}>0$,
$n_{\varepsilon} >0$ and
\( \eta=\eta_{\varepsilon}=(1-2\varepsilon)/2T_{\varepsilon} > 0 \)
 such that the number of disjoint simultaneous short (with respect
to $T$) shadowing time intervals associated to any critical values,
and occurring before time \( n \), is \(
\geq \eta n \) for every \( n\geq
n_{\varepsilon} \).
\end{lemma}

\begin{proof}
Suppose that a point \( x \) satisfies condition
\eqref{TSR}. Then for all \( \varepsilon>0 \) sufficiently
small, there exist \( n_{\varepsilon}, T_{\varepsilon} \) such that
\begin{equation}\label{small}
\sum_{\substack{ 1\leq n_{i}^{\pm} \leq n \\
n_{i+1}^{\pm}-n_{i}^{\pm}>T}}n_{i+1}^{\pm}-n_{i}^{\pm} =
\sum_{\substack{
1\leq n_{i}^{\pm} \leq n \\
s(x^{n_{i}^{\pm}}) > T }} s(x^{n_{i}^{\pm}}) \leq
\sum_{\substack{
1\leq j \leq n \\
s(x^{j}) > T }} s(x^{j})
< \varepsilon n
\end{equation}
for all \( T> T_{\varepsilon} \) and \( n> n_{\varepsilon} \).
Indeed, the equality follows by \eqref{septime}, the first inequality
follows
simply because the summation on the left is over a subset of times
of the summation on the right, and the final inequality follows
directly from \eqref{TSR}.
In particular we have
\[
\sum_{\substack{ 1\leq n_{i}^{\pm} \leq n \\
n_{i+1}^{\pm}-n_{i}^{\pm}\leq T}}n_{i+1}^{\pm}-n_{i}^{\pm} \geq
(1-\varepsilon) n.
\]
The calculation holds for both left
and right shadowing times independently and therefore there exists
at least \( (1-2\varepsilon) n \) iterates which belong to a small
gap for both the left and the right sequence of shadowing times
simultaneously.
Each such iterate is therefore contained in an interval \( \mathcal
T \) of simultaneous short shadowing times and thus it is
possible to define at least \( (1-2\varepsilon) n
/2T_{\varepsilon} \) disjoint such intervals.

Finally, notice that if some point \( x \) satisfies
 condition \eqref{TSR} then any other point on the orbit of \( x \)
 also satisfies \eqref{TSR}.  In particular these conclusions
 hold for all critical values as required.

\end{proof}

\subsection{Exponential shrinking at simultaneous short shadowing times}\label{step1}

\begin{lemma}\label{expshr}
For any integer $T \geq 1$, there exists a constant
$0<\gamma=\gamma (T) <1$ such that for any critical value \( c^{1} \)
and any associated pair \( [n_{i-1}^{-},
n_{i}^{-}] \) and \( [n_{j-1}^{+}, n_{j}^{+}] \) of simultaneous
left and right short gaps we have
$$\frac{|\hat{I}_{-}^{(n_{i}^{-})}|}{|\hat{I}_{-}^{(n_{i-1}^{-})}|}
<\gamma \quad\text{ and }\quad
\frac{|\hat{I}_{+}^{(n_{j}^{+})}|}{|\hat{I}_{+}^{(n_{j-1}^{+})}|}
<\gamma
$$
\end{lemma}

\begin{proof}
    Let \( \mathcal T=[k, \tilde k] \) be the interval containing
  \( [n_{i-1}^{-},
n_{i}^{-}] \) and \( [n_{j-1}^{+}, n_{j}^{+}] \) and satisfying \(
\tilde k - k < T \). We claim first of all that there exists a \(
\delta=\delta(T)>0
    \) such that
    \[
| f^{k}(\hat{I}_{-}^{(n_{i-1}^{-})}\setminus
\hat{I}_{-}^{(n_{i}^{-})})| \geq \delta, \quad
|f^{k}(\hat{I}_{+}^{(n_{j-1}^{+})}\setminus
\hat{I}_{+}^{(n_{j}^{+})})| \geq \delta.
    \]
 To see this, let \( \mathcal C^{(T)}= \{f^{i}(\mathcal C)\}_{i=-T}^{T} \)
 denote the set of all images and preimages of the critical set
 between time \( -T \) and \( +T \) and define \(
 \delta > 0
 \) to be the minimum distance between any two points in this set.
 Now let \( [x,y]= \hat{I}_{-}^{(n_{i-1}^{-})}\setminus \hat{I}_{-}^{(n_{i}^{-})} \)
 and \( [z,w]= \hat{I}_{+}^{(n_{j-1}^{+})}\setminus
    \hat{I}_{+}^{(n_{j}^{+})}\). Then, by construction we have
 \[ \{f^{n_{i-1}^{-}}(x),
    f^{n_{i}^{-}}(y), f^{n_{j-1}^{+}}(w), f^{n_{j}^{+}}(z) \} \subset \mathcal C
 \]
 and therefore, since  \( n_{i-1}^{-}, n_{i}^{-}, n_{j-1}^{+}, n_{j}^{+} \)
 are all within \(
    T \) iterates of \( k \) we have
 \[
 \{f^{k}(x), f^{k}(y), f^{k}(z), f^{k}(w)\}\subset \mathcal C^{(T)}.
 \]
 These are the endpoints of the intervals
\( f^{k}(\hat{I}_{-}^{(n_{i-1}^{-})}\setminus
\hat{I}_{-}^{(n_{i}^{-})}) \) and \(
f^{k}(\hat{I}_{+}^{(n_{j-1}^{+})}\setminus
\hat{I}_{+}^{(n_{j}^{+})}) \)
    and thus the claim follows. (See Figure 4).

\begin{center}
\includegraphics[height=8cm]{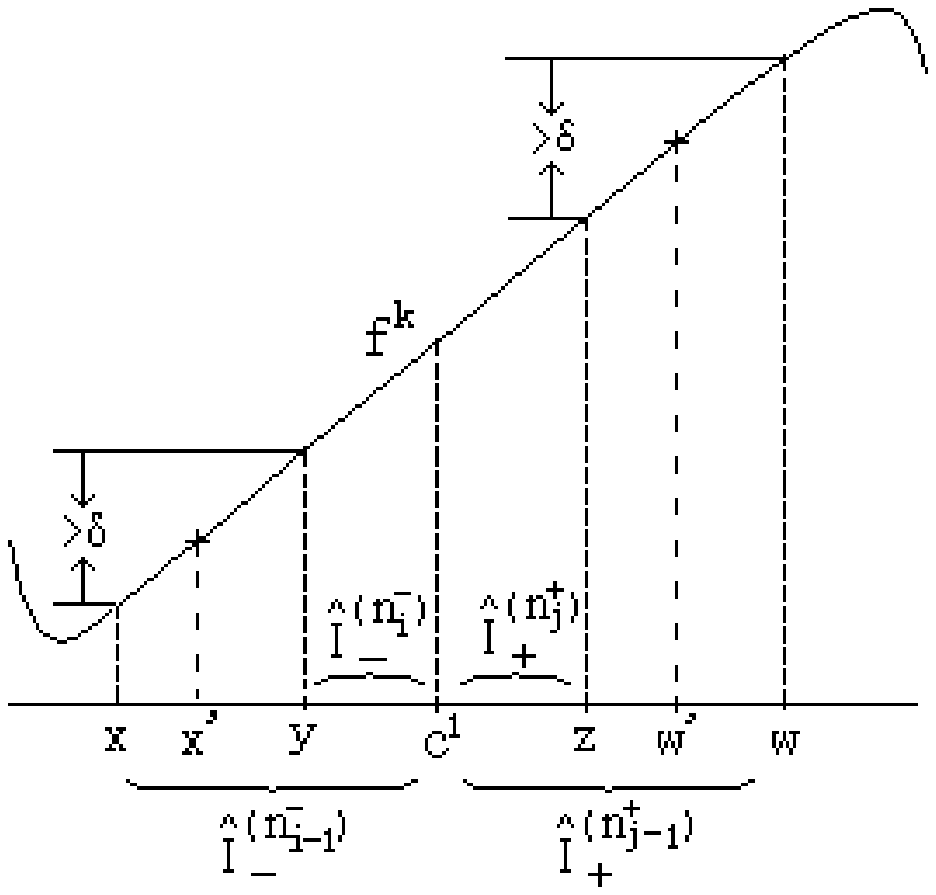}
\makeatletter\def\@captype{figure}\makeatother
\caption{Exponential shrinking at simultaneous short shadowing
times.}
\end{center}

 We now choose two points \( x'\in (x, y) \) and \( w'\in (z,w) \)
    such that \( |f^{k}[x,x']|=|f^{k}[x',y]| \geq \delta/2
    \) and \(|f^{k}[z,w']|=|f^{k}[w',w]| \geq   \delta/2 \).

 By \eqref{monotone} we observe that the map \( f^{k} \) is monotone on
 \( [x, w] \) and therefore the two extreme intervals
 \( [x,x'] \) and \( [w', w] \) form the \emph{Koebe space}
 which guarantees uniformly bounded distortion in \( [x', w'] \)
 (see \cite{MelStr93}). In particular the proportion between the
 lengths of \( [x', y] \) and \( [y, c^1] \) is uniformly comparable
 to the proportion between the length of \( f^{k}[x',y] \) and \(
 f^{k}[y,c^{1}] \) and therefore we have
 \begin{align*}
 \frac{|\hat{I}_{-}^{(n_{i-1}^{-})}|}{|\hat{I}_{-}^{(n_{i}^{-})}|}
 &=\frac{|[x, c^1]|}{|[y, c^1]|} \geq \frac{|[x', c^{1}]| }{|[y, c^{1}]| }
\\ &\geq 1+ \frac{|[x', y]| }{|[y, c^{1}]| }
\geq 1 +
\mathcal \Gamma \frac{ |f^{k}[x',y]| }{|f^{k}[y,c^{1}]| } \geq 1+
\frac{\Gamma\delta}{2|I|}
 \end{align*}
 for a distortion constant \( \Gamma \) which depends only on \( \delta
 \). This proves the first inequality in the
 statement of the Lemma with \( \gamma=(1+\frac{\Gamma\delta}{2|I|} )^{-{1}} \).
 Exactly the same argument gives the second inequality.
 \end{proof}

\subsection{Collet-Eckmann}
\label{global}

We are now ready to show that all critical values admit exponentially
growing derivative along their orbits. More specifically we shall
prove the following
\begin{lemma}\label{CEproof}
Suppose that every critical value satisfies (TSR).
Then there exist constants \( K, \lambda > 0 \) such that for every
critical value and every \( n\geq 1 \) we have
$$r_n :=\min \left\{ \frac{|f^n(\hat I_{-}^{(n)})|}{|\hat I_{-}^{(n)}|},
\frac{|f^n(\hat
I^{(n)}_{+})|}{|\hat I_{+}^{(n)}|} \right\}\geq Ke^{\lambda n}.$$
\end{lemma}

The exponential growth of the derivative then follows from the
so-called \emph{minimum principle} for maps with negative Schwarzian
derivative:

\begin{minprin}\label{negsch}
\cite{MelStr93} Let $f$ be an S-multimodal map and $[a, b]\subset I$
be a subinterval. For $i\geq 1$, if $f^i|_{[a, b]}$ is a
diffeomorphism, then
% $$\frac{|f^i(x)-f^i(y)|}{|x-y|}\geq \min \{\frac{|f^i(b)-f^i(x)|}{|b-x|},
% \; \frac{|f^i(x)-f^i(a)|} {|x-a|}\},$$ especially,
$$|Df^i(x)|\geq \min \left\{ \frac{|f^i(b)-f^i(x)|}{|b-x|},\;
\frac{|f^i(x)-f^i(a)|}{|x-a|}\right\}$$
for any $x \in (a, b)$.
\end{minprin}

We therefore get the following
\begin{corollary}
    Suppose that every critical value satisfies (TSR).
Then there exist constants \( K, \lambda > 0 \) such that for every
critical value \( c^{1} \) and every \( n\geq 1 \) we have
\[
|Df^{i}(c^{1})| \geq K e^{\lambda n}.
\]
\end{corollary}

We begin with a sublemma which follow easily from the estimates
obtained above.

\begin{sublemma}\label{shrink}
     For all \( n\geq 1  \) we have
    \[
    |\hat I_{\pm}^{(n)}| \leq e^{- \eta \gamma n}.
\]
\end{sublemma}
\begin{proof}
Follows from
Lemmas \ref{posden} and  \ref{expshr}.
\end{proof}

\begin{proof}[Proof of Lemma \ref{CEproof}]

We start by considering a (left or right) shadowing
time $n_i^{\pm}$. Let $\hat I^{(n_i^{\pm})}_{\pm}=[x, c^1]$ and recall that
then $f^{n_i^{\pm}}(\hat I_{\pm}^{(n_i^{\pm})})=(\hat c,
c^{n_i^{\pm}+1})$ for some critical point $\hat c$.
 By lemma
\ref{TSRSR} we have
$$d(c^j)\geq e^{-s(c^j)/\underline \kappa}.$$
Therefore
\begin{equation}\label{shaexp}
|f^{n_i^{\pm}}(\hat I_{\pm}^{(n_i^{\pm})})| =
|c^{n_i^{\pm}+1}-\hat
c| \geq e^{-s(c^{n_i^{\pm}+1})/\underline\kappa }
\end{equation}
for every  critical value and every shadowing time $n_i^{\pm}\geq
1$.

We now let $n_i^{\pm}$ be the smallest (left or right) shadowing
time larger than $n$. Then $n_{i-1}^{\pm}\leq n <n_i^{\pm}$ and so
$\hat I_{\pm}^{(n_i^{\pm})} \subset \hat I^{(n)}_{\pm}=\hat
I_{\pm}^{(n_{i-1}^{\pm})}$. By \eqref{monotone} we know that
$f^{n_i^{\pm}}|_{\hat I_{\pm}^{(n_{i-1}^{\pm})}}$ is monotone, so
we have by \eqref{shaexp},
\begin{equation}\label{NNI}
\begin{split}
|f^n(\hat I_{\pm}^{(n)})| & \geq
D^{-(n_i^{\pm}-n)}|f^{n^{\pm}_{i}}(\hat I_{\pm}^{(n)})|=
D^{-(n_i^{\pm}-n)}|f^{n^{\pm}_{i}}(\hat
I_{\pm}^{(n_{i-1}^{\pm})})| \\
&\geq D^{-(n_i^{\pm}-n)}|f^{n^{\pm}_{i}}(\hat
I_{\pm}^{(n_{i}^{\pm})})|.
\end{split}
\end{equation}
Then, by \eqref{shaexp} and sublemma \ref{shrink}, we obtain
\begin{equation}
\begin{split}
|f^n(\hat I_{\pm}^{(n)})| \geq
D^{-(n_i^{\pm}-n)}e^{-s(c^{n_i^{\pm}+1})/\underline\kappa } & \geq
   D^{-(n_i^{\pm}-n)}e^{-s(c^{n_i^{\pm}+1})/\underline \kappa }
        e^{\eta \gamma n }|\hat I_{\pm}^{(n)}|
        \\
         & =D^{-(n_i^{\pm}-n)}
e^{-s(c^{n_i^{\pm}+1})/\underline\kappa+ \eta \gamma n }
 |\hat
I_{\pm}^{(n)}|.
 \end{split}
\end{equation}
It is therefore sufficient to show that we can choose constants \(
K, \lambda > 0 \) so that
\[
r_{n}\geq \ D^{-(n_i^{\pm}-n)}
e^{-s(c^{n_i^{\pm}+1})/\underline\kappa+ \eta \gamma n }   \geq
Ke^{\lambda n}.
\]
Since $\lim_{n\rightarrow +\infty}{s(c^n)}/{n}=0$, we can choose
$N$ sufficiently large such that
$s(c^{n_i^{\pm}+1})/\underline\kappa  < \eta \gamma
(n_i^{\pm}+1)/2$ when $n_i^{\pm} \geq N$. So
$$ r_{n}\geq \ D^{-(n_i^{\pm}-n)}e^{\eta \gamma n-\eta \gamma
(n_i^{\pm}+1)/2}.$$ We start by taking $\tilde \lambda
=\eta\gamma/2 $, fixing some arbitrary \( \tilde \lambda
> \lambda > 0 \) and writing \( \tilde \lambda=\lambda+
(\tilde \lambda-\lambda) \) and \( n_{i}^{\pm}=(n_{i}^{\pm}-n)+n
\). Then we have \( e^{\tilde \lambda n_{i}^{\pm}} = e^{\lambda
n}e^{(\tilde \lambda - \lambda)n}e^{-\tilde \lambda
(n-n_{i}^{\pm})} \) and therefore
\[r_{n}\geq \ D^{-(n_i^{\pm}-n)}e^{\eta \gamma n-\eta \gamma
(n_i^{\pm}+1)/2}=D^{-(n_i^{\pm}-n)}e^{-\tilde \lambda
(n_i^{\pm}-n)}e^{-\tilde \lambda}e^{(\tilde \lambda-
\lambda)n}e^{\lambda n}.
\]
It remains to show that
\begin{equation}
\label{bgtnc} (De^{\tilde\lambda})^{-(n_{i}^{\pm}-n)}
e^{-\tilde\lambda} e^{(\tilde\lambda - \lambda) n} \geq K
\end{equation}
This follows by the crucial facts that \(
n_{i}^{\pm}-n_{i-1}^{\pm}> n-n_{i-1}^{\pm} \) and
\begin{equation}\label{crucial}
\lim_{n_{i-1}^{\pm}\rightarrow +\infty
}\frac{s(c^{n_{i-1}^{\pm}})}{n_{i-1}^{\pm}}=
\lim_{n_{i-1}^{\pm}\rightarrow +\infty
}\frac{n_{i}^{\pm}-n_{i-1}^{\pm}}{n_{i-1}^{\pm}} =0
\end{equation}
by TSR and \eqref{septime}. Indeed, We use \eqref{crucial} which
implies that the left side of \eqref{bgtnc} is bigger than $1$ for
all \( n_{i}^{\pm}\geq \tilde N \) for some sufficiently large \(
\tilde N \) depending only on the map and the constants \( \lambda
\) and \( \tilde\lambda \) and thus, again using \eqref{crucial},
for all \( n \geq N \)  for some sufficiently large \( N \) also
depending only on the same quantities as \( \tilde N \). To take
care of smaller values of \( n \) it is sufficient to choose \( K
\) sufficiently small depending only on the value of \( r_{n} \)
for \( n\leq N \). This completes the proof.

% First we take $N_1 \geq N_0$ such
% that $\tilde K \tilde \lambda^{N_1}
% >1$, and choose $\lambda$ such that $0< \lambda < \tilde \lambda$
% and $\lambda ^{N_1} < \tilde \lambda^{N_1}$. As
% $\lim_{n\rightarrow +\infty }s(c^n)/n=0$, we choose $N_2>N_1$
% sufficient large such that for any shadowing time
% $n_i^{\pm}>N_2$, we have
% \begin{equation}\label{LARN}
% e^{(\tilde \lambda-\lambda)n_i^{\pm}}\geq D^{n_{i+1}^{\pm}-n_i^{\pm}}.
% \end{equation}
% Finally we take $K=\min \{\tilde K, \min_{1 \leq j \leq N_2} r_j
% \}/e^{N_2\lambda}$. By choosing the above constants, we have
% proved $r_i \geq Ke^{\lambda i}$ for $1 \leq i \leq N_2$.
% Using \eqref{LARN} we
% obtain $|f^n(\hat I_{-}^{(n)})| \geq Ke^{\lambda n} |\hat
% I_{-}^{(n)}|$ for any $n>N_2$. Similar argument shows the other
% inequality, we have $r_n \geq K e^{\lambda n}$ and complete the
% proof.

\end{proof}

\section*{Acknowledgements}
Lanyu Wang would like to thank his PhD supervisor Professor Zhifen
Zhang for her advise and encouragement, and Professors Zhujun Jing
and Lan Wen for their help and encouragement for many years.

\end{document}